# Support Vector Machines with Applications[1]

**Javier M. Moguerza and Alberto Muñoz**

*Abstract.* Support vector machines (SVMs) appeared in the early nineties as optimal margin classifiers in the context of Vapnik's statistical learning theory. Since then SVMs have been successfully applied to real-world data analysis problems, often providing improved results compared with other techniques. The SVMs operate within the framework of regularization theory by minimizing an empirical risk in a well-posed and consistent way. A clear advantage of the support vector approach is that sparse solutions to classification and regression problems are usually obtained: only a few samples are involved in the determination of the classification or regression functions. This fact facilitates the application of SVMs to problems that involve a large amount of data, such as text processing and bioinformatics tasks. This paper is intended as an introduction to SVMs and their applications, emphasizing their key features. In addition, some algorithmic extensions and illustrative real-world applications of SVMs are shown.

## 1. INTRODUCTION

In the last decade, support vector machines (SVMs) have increasingly turned into a standard methodology in the computer science and engineering communities. As Breiman [12] pointed out, these communities are often involved in the solution of con-sulting and industrial data analysis problems. The usual starting point is a sample data set $\{(\mathbf{x}_i, \mathbf{y}_i) \in X \times Y\}_{i=1}^n$, and the goal is to "learn" the relationship between the $\mathbf{x}$ and $\mathbf{y}$ variables. The variable $X$ may be, for instance, the space of $20 \times 20$ binary matrices that represent alphabetic uppercase characters and $Y$ would be the label set $\{1, \ldots, 27\}$. Similarly, $X$ may be $\mathbb{R}^{10,000}$, the space corresponding to a document data base with a vocabulary of 10,000 different words. In this case $Y$ would be the set made up of a finite number of predefined semantic document classes, such as statistics, computer science, sociology and so forth. The main goal in this context usually is predictive accuracy, and in most cases it is not possible to assume a parametric form for the probability distribution $p(\mathbf{x}, \mathbf{y})$. Within this setting many practitioners concerned with providing practical solutions to industrial data analysis problems put more emphasis on algorithmic modeling than on data models. However, a solely algorithmic point of view can lead to procedures with a black box behavior, or even worse, with a poor response to the bias–variance dilemma. Neural networks constitute a paradigmatic example of this approach. The (semi-parametric) model implemented by neural networks is powerful enough to approximate continuous functions with arbitrary precision. On the other hand, neural network parameters are very hard to tune and interpret, and statistical inference is usually not possible [51].

*Javier M. Moguerza is Associate Professor, School of Engineering, University Rey Juan Carlos, c/ Tulipan s/n, 28933 Mostoles, Spain e-mail:*
*javier.moguerza@urjc.es. Alberto Muñoz is Associate Professor, Department of Statistics, University Carlos III, c/Madrid 126, 28903 Getafe, Spain e-mail:*
*alberto.munoz@uc3m.es.*
[1]*Discussed in* 10.1214/08834230600000457,
10.1214/088342306000000466, 10.1214/088342306000000475
and 10.1214/088342306000000484; *rejoinder*
10.1214/088342306000000501.

The SVMs provide a compromise between the parametric and the pure nonparametric approaches: As in linear classifiers, SVMs estimate a linear decision function, with the particularity that a previous mapping of the data into a higher-dimensional feature space may be needed. This mapping is characterized by the choice of a class of functions known as kernels. The support vector method was introduced by Boser, Guyon and Vapnik [10] at the Computational Learning Theory (COLT92) ACM Conference. Their proposal subsumed into an elegant and theoretically well founded algorithm two seminal ideas, which had already individually appeared throughout previous years: the use of kernels and their geometrical interpretation, as introduced by Aizerman, Braverman and Rozonoer [1], and the idea of constructing an optimal separating hyperplane in a nonparametric context, developed by Vapnik and Chervonenkis [78] and by Cover [16]. The name "support vector" was explicitly used for the first time by Cortes and Vapnik [15]. In recent years, several books and tutorials on SVMs have appeared. A reference with many historical annotations is the book by Cristianini and Shawe-Taylor [20]. For a review of SVMs from a purely geometrical point of view, the paper by Bennett and Campbell [9] is advisable. An exposition of kernel methods with a Bayesian taste can be read in the book by Herbrich [30]. Concerning the statistical literature, the book by Hastie, Tibshirani and Friedman [28] includes a chapter dedicated to SVMs.

We illustrate the basic ideas of SVMs for the two-group classification problem. This is the typical version and the one that best summarizes the ideas that underlie SVMs. The issue of discriminating more than two groups can be consulted, for instance, in [37].

Consider a classification problem where the discriminant function is nonlinear, as illustrated in Figure 1(a). Suppose we have a mapping $\Phi$ into a "feature space" such that the data under consideration have become linearly separable as illustrated in Figure 1(b). From the infinite number of existing separating hyperplanes, the support vector machine looks for the plane that lies furthermost from both classes, known as the optimal (maximal) margin hyperplane. To be more specific, denote the available mapped sample by $\{(\Phi(\mathbf{x}_i), y_i)\}_{i=1}^n$, where $y_i \in \{-1, +1\}$ indicates the two possible classes. Denote by $\mathbf{w}^T\Phi(\mathbf{x}) + b = 0$ any separating hyperplane in the space of the mapped data equidistant to the nearest point in each class. Under the assumption of separability, we can rescale $\mathbf{w}$ and $b$ so that $|\mathbf{w}^T\Phi(\mathbf{x}) + b| = 1$ for those points in each class nearest to the hyperplane. Therefore, it holds that for every $i \in \{1, \ldots, n\}$,

$$(1.1) \quad \mathbf{w}^T\Phi(\mathbf{x}_i) + b \begin{cases} \geq 1, & \text{if } y_i = +1 \\ \leq -1, & \text{if } y_i = -1. \end{cases}$$

After the rescaling, the distance from the nearest point in each class to the hyperplane is $1/\|\mathbf{w}\|$. Hence, the distance between the two groups is $2/\|\mathbf{w}\|$, which is called the margin. To maximize the margin, the following optimization problem has to be solved:

$$(1.2) \quad \begin{aligned} & \min_{\mathbf{w}, b} & & \|\mathbf{w}\|^2 \\ & \text{subject to (s.t.)} & & \\ & & & y_i(\mathbf{w}^T\Phi(\mathbf{x}_i) + b) \geq 1, \\ & & & i = 1, \ldots, n, \end{aligned}$$

where the square in the norm of $\mathbf{w}$ has been introduced to make the problem quadratic. Notice that, given its convexity, this optimization problem has no local minima. Consider the solution of problem (1.2), and denote it by $\mathbf{w}^*$ and $b^*$. This solution determines the hyperplane in the feature space $D^*(\mathbf{x}) = (\mathbf{w}^*)^T\Phi(\mathbf{x}) + b^* = 0$. Points $\Phi(\mathbf{x}_i)$ that satisfy the equalities $y_i((\mathbf{w}^*)^T\Phi(\mathbf{x}_i) + b^*) = 1$ are called support vectors [in Figure 1(b) the support vectors are the black points]. As we will make clear later, the support vectors can be automatically determined from the solution of the optimization problem. Usually the support vectors represent a small fraction of the sample, and the solution is said to be sparse. The hyperplane $D^*(\mathbf{x}) = 0$ is completely determined by the subsample made up of the support vectors. This fact implies that, for many applications, the evaluation of the decision function $D^*(\mathbf{x})$ is computationally efficient, allowing the use of SVMs on large data sets in real-time environments.

The SVMs are especially useful within ill-posed contexts. A discussion of ill-posed problems from a statistical point of view may be seen in [55]. A common ill-posed situation arises when dealing with data sets with a low ratio of sample size to dimension. This kind of difficulty often comes up in problems such as automatic classification of web pages or microarrays. Consider, for instance, the following classification problem, where the data set is a text data base that contains 690 documents. These documents have been retrieved from the LISA (Library Science Abstracts) and the INSPEC (bibliographic references for physics, computing and engineering research, from the IEE Institute) data bases, using, respectively, the search keywords "library science" (296 records) and "pattern recognition" (394 records). We have selected as data points the terms that occur in at least ten documents, obtaining 982 terms. Hence, the data set is given by a $982 \times 690$ matrix, say $T$, where $T_{ij} = 1$ if term $i$ occurs in document $j$ and $T_{ij} = 0$ otherwise. For each term, we check the number of library science and pattern recognition documents that contain it. The highest value determines the class of the term. This procedure is standard in the field of automatic thesaurus generation (see [5]). The task is to check the performance of the SVM classifier in recovering the class labels obtained by the previous procedure. Notice that we are dealing with about 1000 points in nearly 700 dimensions. We have divided the data set into a training set (80% of the data points) and a test set (20% of the data points). Since the sample is relatively small with respect to the space dimension, it should be easy for any method to find a criterion that separates the training set into two classes, but this does not necessarily imply the ability to correctly classify the test data.

The results obtained using Fisher linear discriminant analysis (FLDA), the $k$-nearest neighbor classifier ($k$-NN) with $k = 1$ and the linear SVM [i.e., taking $\Phi$ as the identity map $\Phi(x) = x$] are shown in Table 1.

TABLE 1
*Classification percentage errors for a two-class text data base*

| Method | Training error | Test error |
|---|---|---|
| FLDA | 0.0% | 31.4% |
| $k$-NN ($k = 1$) | 0.0% | 14.0% |
| Linear SVM | 0.0% | 3.0% |

It is apparent that the three methods have been able to find a criterion that perfectly separates the training data set into two classes, but only the linear SVM shows good performance when classifying new data points. The best result for the $k$-NN method (shown in the table) is obtained for $k = 1$, an unsurprising result, due to the "curse of dimensionality" phenomenon, given the high dimension of the data space. Regarding FLDA, the estimation of the mean vectors and covariance matrices of the groups is problematic given the high dimension and the small number of data points. The SVMs also calculate a linear hyperplane, but are looking for something different—margin maximization, which will only depend on the support vectors. In addition, there is no loss of information caused by projections of the data points. The successful behavior of the support vector method is not casual, since, as we will see below, SVMs are supported by regularization theory, which is particularly useful for the solution of ill-posed problems like the present one.

In summary, we have just described the basics of a classification algorithm which has the following features:

• Reduction of the classification problem to the computation of a linear decision function.

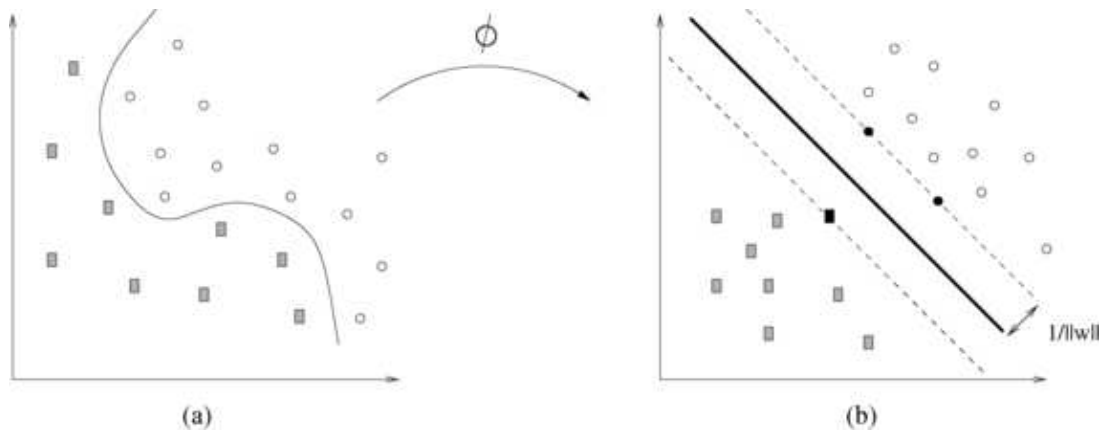

FIG. 1. (a) *Original data in the input space.* (b) *Mapped data in the feature space.*

- Absence of local minima in the SVM optimization problem.
- A computationally efficient decision function (sparse solution).

In addition, in the next sections we will also discuss other important features such as the use of kernels as a primary source of information or the tuning of a very reduced set of parameters.

The rest of the paper is organized as follows. Section 2 shows the role of kernels within the SVM approach. In Section 3 SVMs are developed from the regularization theory perspective and some illustrative examples are given. Section 4 reviews a number of successful SVM applications to real-world problems. In Section 5 algorithmic extensions of SVMs are presented. Finally, in Section 6 some open questions and final remarks are presented.

## 2. THE KERNEL MAPPING

In this section we face one of the key issues of SVMs: how to use $\Phi(\mathbf{x})$ to map the data into a higher-dimensional space. This procedure is justified by Cover's theorem [16], which guarantees that any data set becomes arbitrarily separable as the data dimension grows. Of course, finding such nonlinear transformations is far from trivial. To achieve this task, a class of functions called kernels is used. Roughly speaking, a kernel $K(\mathbf{x}, \mathbf{y})$ is a real-valued function $K : X \times X \to \mathbb{R}$ for which there exists a function $\Phi : X \to Z$, where $Z$ is a real vector space, with the property $K(\mathbf{x}, \mathbf{y}) = \Phi(\mathbf{x})^T \Phi(\mathbf{y})$. This function $\Phi$ is precisely the mapping in Figure 1. The kernel $K(\mathbf{x}, \mathbf{y})$ acts as a dot product in the space $Z$. In the SVM literature $X$ and $Z$ are called, respectively, input space and feature space (see Figure 1).

As an example of such a $K$, consider two data points $\mathbf{x}_1$ and $\mathbf{x}_2$, with $\mathbf{x}_i = (x_{i1}, x_{i2})^T \in \mathbb{R}^2$, and $K(\mathbf{x}_1, \mathbf{x}_2) = (1 + \mathbf{x}_1^T \mathbf{x}_2)^2 = (1 + x_{11}x_{21} + x_{12}x_{22})^2 = \Phi(\mathbf{x}_1)^T \Phi(\mathbf{x}_2)$, where $\Phi(\mathbf{x}_i) = (1, \sqrt{2}x_{i1}, \sqrt{2}x_{i2}, x_{i1}^2, x_{i2}^2, \sqrt{2}x_{i1}x_{i2})$. Thus, in this example $\Phi : \mathbb{R}^2 \to \mathbb{R}^6$. As we will show later, explicit knowledge of both the mapping $\Phi$ and the vector $\mathbf{w}$ will not be needed: we need only $K$ in its closed form.

To be more specific, a kernel $K$ is a positive definite function that admits an expansion of the form $K(\mathbf{x}, \mathbf{y}) = \sum_{i=1}^{\infty} \lambda_i \Phi_i(\mathbf{x}) \Phi_i(\mathbf{y})$, where $\lambda_i \in \mathbb{R}^+$. Sufficient conditions for the existence of such an expansion are given in Mercer's theorem [43]. The function $K(\mathbf{x}, \mathbf{y})$, known as a Mercer's kernel, implicitly defines the mapping $\Phi$ by letting $\Phi(\mathbf{x}) = (\sqrt{\lambda_1}\Phi_1(\mathbf{x}), \sqrt{\lambda_2}\Phi_2(\mathbf{x}), \ldots)^T$.

Examples of Mercer's kernels are the linear kernel $K(\mathbf{x}, \mathbf{y}) = \mathbf{x}^T \mathbf{y}$, polynomial kernels $K(\mathbf{x}, \mathbf{y}) = (c + \mathbf{x}^T \mathbf{y})^d$ and the Gaussian kernel $K_c(\mathbf{x}, \mathbf{y}) = e^{-\|\mathbf{x}-\mathbf{y}\|^2/c}$. In the first case, the mapping is the identity. Polynomial kernels map the data into finite-dimensional vector spaces. With the Gaussian kernel, the data are mapped onto an infinite dimensional space $Z = \mathbb{R}^\infty$ (all the $\lambda_i \neq 0$ in the kernel expansion; see [63] for the details).

Given a kernel $K$, we can consider the set of functions spanned by finite linear combinations of the form $f(\mathbf{x}) = \sum_j \alpha_j K(\mathbf{x}_j, \mathbf{x})$, where the $\mathbf{x}_j \in X$. The completion of this vector space is a Hilbert space known as a reproducing kernel Hilbert space (RKHS) [3]. Since $K(\mathbf{x}_j, \mathbf{x}) = \Phi(\mathbf{x}_j)^T \Phi(\mathbf{x})$, the functions $f(\mathbf{x})$ that belong to a RKHS can be expressed as $f(\mathbf{x}) = \mathbf{w}^T \Phi(\mathbf{x})$, with $\mathbf{w} = \sum_j \alpha_j \Phi(\mathbf{x}_j)$, that is, $f(\mathbf{x}) = 0$ describes a hyperplane in the feature space determined by $\Phi$ [as the one illustrated in Figure 1(b)]. Thus, reproducing kernel Hilbert spaces provide a natural context for the study of hyperplanes in feature spaces through the use of kernels like those introduced in Section 1. Without loss of generality, a constant $b$ can be added to $f$ (see [64] for a complete discussion), taking the form

$$(2.1) \qquad f(\mathbf{x}) = \sum_j \alpha_j K(\mathbf{x}_j, \mathbf{x}) + b.$$

Equation (2.1) answers the question of how to use $\Phi(\mathbf{x})$ to map the data onto a higher-dimensional space: Since $f(\mathbf{x})$ can be evaluated using expression (2.1) [in which only the kernel values $K(\mathbf{x}_j, \mathbf{x})$ are involved], $\Phi$ acts implicitly through the closed form of $K$. In this way, the kernel function $K$ is employed to avoid an explicit evaluation of $\Phi$ (often a high-dimensional mapping). This is the reason why knowledge of the explicit mapping $\Phi$ is not needed.

As we will show in the next section, SVMs work by minimizing a regularization functional that involves an empirical risk plus some type of penalization term. The solution to this problem is a function that has the form (2.1). This optimization process necessarily takes place within the RKHS associated with the kernel $K$. The key point in this computation is the way in which SVMs select the weights $\alpha_j$ in (2.1) (the points $\mathbf{x}_j$ are trivially chosen as the sample data points $\mathbf{x}_i$). A nice fact is that the estimation of these weights, which determine the decision function in the RKHS, is reduced to the solution of a smooth and convex optimization problem.

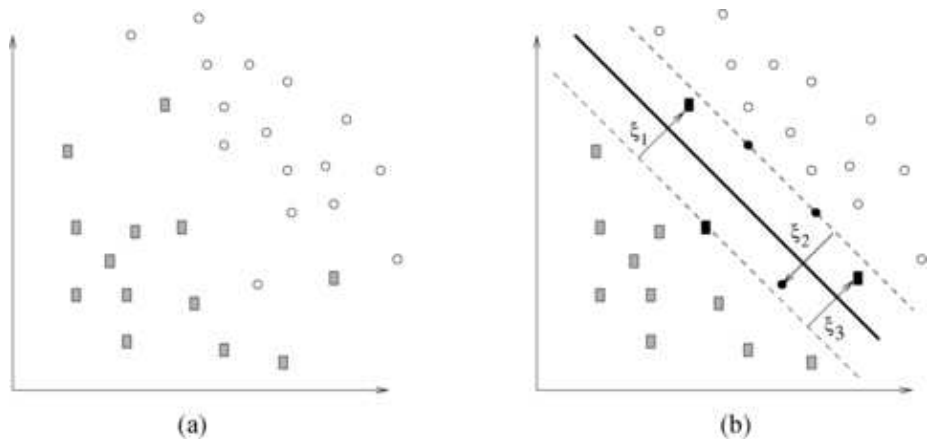

Fig. 2. (a) *Nonseparable mapped data in the feature space.* (b) *Normalized hyperplane for the data in* (a).

## 3. SUPPORT VECTOR MACHINES: A REGULARIZATION METHOD

In Section 1 we introduced the formulation of SVMs for the situation illustrated in Figure 1(b), where the mapped data have become linearly separable. We consider now the more general case where the mapped data remain nonseparable. This situation is illustrated in Figure 2(a). The SVMs address this problem by finding a function $f$ that minimizes an empirical error of the form $\sum_{i=1}^{n} L(y_i, f(\mathbf{x}_i))$, where $L$ is a particular loss function and $(\mathbf{x}_i, y_i)_{i=1}^{n}$ is the available data sample. There may be an infinite number of solutions, in which case the problem is ill-posed. Our aim is to show how SVMs make the problem well-posed. As a consequence, the decision function calculated by the SVM will be unique, and the solution will depend continuously on the data.

The specific loss function $L$ used within the SVM approach is $L(y_i, f(\mathbf{x}_i)) = (1 - y_i f(\mathbf{x}_i))_+$, with $(x)_+ = \max(x, 0)$. This loss function is called hinge loss and is represented in Figure 3. It is zero for well classified points with $|f(\mathbf{x}_i)| \geq 1$ and is linear otherwise. Hence, the hinge loss function does not penalize large values of $f(\mathbf{x}_i)$ with the same sign as $y_i$ (understanding large to mean $|f(\mathbf{x}_i)| \geq 1$).

This behavior agrees with the fact that in classification problems only an estimate of the classification boundary is needed. As a consequence, we only take into account points such that $L(y_i, f(\mathbf{x}_i)) > 0$ to determine the decision function.

To reach well-posedness, SVMs make use of regularization theory, for which several similar approaches have been proposed [33, 60, 73]. The widest used setting minimizes Tikhonov's regularization functional [73], which consists of solving the optimization problem

$$(3.1) \quad \min_{f \in H_K} \frac{1}{n} \sum_{i=1}^{n} (1 - y_i f(\mathbf{x}_i))_+ + \mu \|f\|_K^2,$$

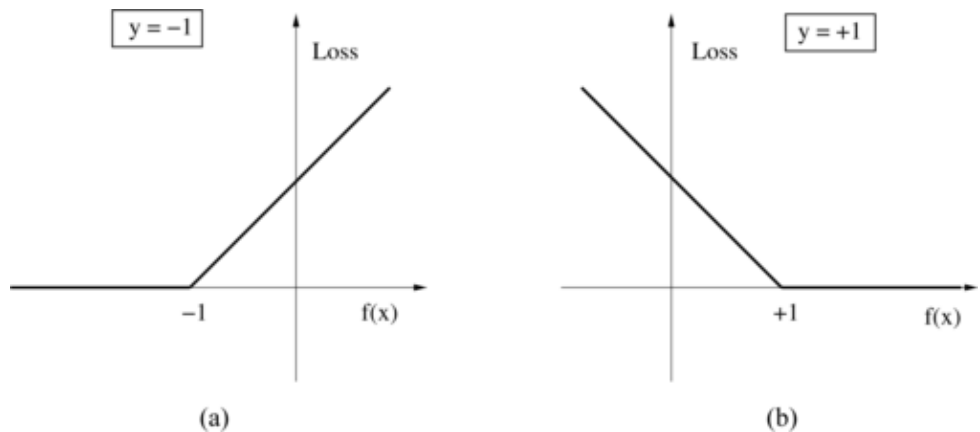

Fig. 3. *Hinge loss function* $L(y_i, f(\mathbf{x}_i)) = (1 - y_i f(\mathbf{x}_i))_+$: (a) $L(-1, f(\mathbf{x}_i))$; (b) $L(+1, f(\mathbf{x}_i))$.

where $\mu > 0$, $H_K$ is the RKHS associated with the kernel $K$, $\|f\|_K$ denotes the norm of $f$ in the RKHS and $\mathbf{x}_i$ are the sample data points. Given that $f$ belongs to $H_K$, it takes the form $f(\cdot) = \sum_j \alpha_j K(\mathbf{x}_j, \cdot)$. As in Section 2, $f(\mathbf{x}) = 0$ is a hyperplane in the feature space. Using the reproducing property $\langle K(\mathbf{x}_j, \cdot), K(\mathbf{x}_l, \cdot) \rangle_K = K(\mathbf{x}_j, \mathbf{x}_l)$ (see [3]), it holds that $\|f\|_K^2 = \langle f, f \rangle_K = \sum_j \sum_l \alpha_j \alpha_l K(\mathbf{x}_j, \mathbf{x}_l)$.

In (3.1) the scalar $\mu$ controls the trade-off between the fit of the solution $f$ to the data (measured by $L$) and the approximation capacity of the function space that $f$ belongs to (measured by $\|f\|_K$). It can be shown [11, 48] that the space where the solution is sought takes the form $\{f \in H_K : \|f\|_K^2 \le (\sup_{y \in Y} L(y, 0))/\mu\}$, a compact ball in the RKHS. Note that the larger $\mu$ is, the smaller is the ball and the more restricted is the search space. This is the way in which regularization theory imposes compactness in the RKHS. Cucker and Smale [21] showed that imposing compactness on the space assures well-posedness of the problem and, thus, uniqueness of the solution (refer to the Appendix for details).

The solution to problem (3.1) has the form $f(\mathbf{x}) = \sum_{i=1}^{n} \alpha_i K(\mathbf{x}_i, \mathbf{x}) + b$, where $\mathbf{x}_i$ are the sample data points, a particular case of (2.1). This result is known as the representer theorem. For details, proofs and generalizations, refer to [36, 67] or [18]. It is immediate to show that $\|f\|_K^2 = \|\mathbf{w}\|^2$, where $\mathbf{w} = \sum_i^n \alpha_i \Phi(\mathbf{x}_i)$. Given this last result, problem (3.1) can be restated as

$$(3.2) \quad \min_{\mathbf{w}, b} \frac{1}{n} \sum_{i=1}^{n} (1 - y_i(\mathbf{w}^T \Phi(\mathbf{x}_i) + b))_+ + \mu \|\mathbf{w}\|^2.$$

It is worth mentioning that the second term in (3.2) coincides with the term in the objective function of (1.2). Problems (3.1) and (3.2) review some of the key issues of SVMs enumerated at the end of Section 1: Through the use of kernels, the a priori problem of estimating a nonlinear decision function in the input space is transformed into the a posteriori problem of estimating the weights of a hyperplane in the feature space.

Because of the hinge loss function, problem (3.2) is nondifferentiable. This lack of differentiability implies a difficulty for efficient optimization techniques; see [7] or [47]. Problem (3.2) can be turned smooth by straightforwardly formulating it as (see [41])

$$\min_{\mathbf{w}, b, \xi} \quad \frac{1}{2} \|\mathbf{w}\|^2 + C \sum_{i=1}^{n} \xi_i$$

$$(3.3) \quad \begin{aligned} \text{s.t.} \quad & y_i(\mathbf{w}^T \Phi(\mathbf{x}_i) + b) \ge 1 - \xi_i, \\ & i = 1, \dots, n, \\ & \xi_i \ge 0, \quad i = 1, \dots, n, \end{aligned}$$

where $\xi_i$ are slack variables introduced to avoid the nondifferentiability of the hinge loss function and $C = 1/(2\mu n)$. This is the most widely used SVM formulation.

The slack variables $\xi_i$ allow violations of constraints (1.1), extending problem (1.2) to the nonseparable case [problem (1.2) would not be solvable for nonseparable data]. The slack variables guarantee the existence of a solution. The situation is shown in Figure 2(b), which constitutes a generalization of Figure 1(b). Notice that problem (1.2) is a particular case of problem (3.3). To be more specific, if the mapped data become separable, problem (1.2) is equivalent to problem (3.3) when, at the solution, $\xi_i = 0$. Intuitively, we want to solve problem (1.2) and, at the same time, minimize the number of nonseparable samples, that is, $\sum_i \#(\xi_i > 0)$. Since the inclusion of this term would provide a nondifferentiable combinatorial problem, the smooth term $\sum_{i=1}^{n} \xi_i$ appears instead.

We have deduced the standard SVM formulation (3.3) via the use of regularization theory. This framework guarantees that the empirical error for SVMs converges to the expected error as $n \to \infty$ [21], that is, the decision functions obtained by SVMs are statistically consistent. Therefore, the separating hyperplanes obtained by SVMs are neither arbitrary nor unstable. This remark is pertinent since Cover's theorem (which guarantees that any data set becomes arbitrarily separable as the data dimension grows) could induce some people to think that SVM classifiers are arbitrary.

By standard optimization theory, it can be shown that problem (3.3) is equivalent to solving

$$\min_{\lambda} \quad \frac{1}{2} \sum_{i=1}^{n} \sum_{j=1}^{n} \lambda_i \lambda_j y_i y_j K(\mathbf{x}_i, \mathbf{x}_j) - \sum_{i=1}^{n} \lambda_i$$

$$(3.4) \quad \begin{aligned} \text{s.t.} \quad & \sum_{i=1}^{n} y_i \lambda_i = 0, \\ & 0 \le \lambda_i \le C, \quad i = 1, \dots, n. \end{aligned}$$

The $\lambda_i$ variables are the Lagrange multipliers associated with the constraints in (3.3). This problem is known in optimization theory as the dual problem of (3.3) [7]. This problem is convex and quadratic and, therefore, every local minimum is a global minimum.

In practice, this is the problem to solve, and efficient methods specific for SVMs have been developed (see [34, 58, 61]).

Let the vector $\lambda^*$ denote the solution to problem (3.4). Points that satisfy $\lambda_i^* > 0$ are the support vectors (shown in black in Figure 2(b) for the nonseparable case). It can be shown that the solution to problem (3.3) is $\mathbf{w}^* = \sum_{i=1}^n \lambda_i^* y_i \Phi(\mathbf{x}_i)$ and

$$(3.5) \quad \begin{aligned} b^* = & -\frac{\sum_{i=1}^n \lambda_i^* y_i K(\mathbf{x}_i, \mathbf{x}^+)}{2} \\ & + \frac{\sum_{i=1}^n \lambda_i^* y_i K(\mathbf{x}_i, \mathbf{x}^-)}{2}, \end{aligned}$$

where $\mathbf{x}^+$ and $\mathbf{x}^-$ are, respectively, two support vectors in classes $+1$ and $-1$ such that their associated Lagrange multipliers $\lambda^+$ and $\lambda^-$ hold so that $0 < \lambda^+ < C$ and $0 < \lambda^- < C$.

The desired decision function, which determines the hyperplane $(\mathbf{w}^*)^T \Phi(\mathbf{x}) + b^* = 0$, takes the form

$$(3.6) \quad \begin{aligned} D^*(\mathbf{x}) &= (\mathbf{w}^*)^T \Phi(\mathbf{x}) + b^* \\ &= \sum_{i=1}^n \lambda_i^* y_i K(\mathbf{x}_i, \mathbf{x}) + b^*. \end{aligned}$$

Equations (3.5) and (3.6) show that $D^*(\mathbf{x})$ is completely determined by the subsample made up by the support vectors, the only points in the sample for which $\lambda_i^* \neq 0$. This definition of support vector is coherent with the geometrical one given in Section 1. The reason is that Lagrange multipliers $\lambda_i^*$ must fulfill the strict complementarity conditions (see [7]), that is, $\lambda_i^*(D^*(\mathbf{x}_i) - 1 + \xi_i) = 0$, where either $\lambda_i^* = 0$ or $D^*(\mathbf{x}_i) = 1 - \xi_i$. Therefore, if $\lambda_i^* \neq 0$, then $D^*(\mathbf{x}_i) = 1 - \xi_i$ and $\mathbf{x}_i$ is one of the points that defines the decision hyperplane [one of the black points in Figure 2(b)]. Often the support vectors are a small fraction of the data sample and, as already mentioned, the solution is said to be sparse. This property is due to the use of the hinge loss function.

Note that problem (3.4) and equation (3.6) depend only on kernel evaluations of the form $K(\mathbf{x}, \mathbf{y})$. Therefore, the explicit mapping $\Phi$ is not needed to solve the SVM problem (3.4) or to evaluate the decision hyperplane (3.6). In particular, even when the kernel corresponds to an infinite-dimensional space (for instance, the Gaussian kernel), there is no problem with the evaluation of $\mathbf{w}^* = \sum_{i=1}^n \lambda_i^* y_i \Phi(\mathbf{x}_i)$, which is not explicitly needed. In practice, $D^*(\mathbf{x})$ is evaluated using the right-hand side of equation (3.6).

## 3.1 SVMs and the Optimal Bayes Rule

The results in the previous section are coherent with the ones obtained by Lin [40], which state that the support vector machine classifier approaches the optimal Bayes rule and its generalization error converges to the optimal Bayes risk.

Consider a two-group classification problem with classes $+1$ and $-1$ and, to simplify, assume equal costs of misclassification. Under this assumption, the expected misclassification rate and the expected cost coincide. Let $p_1(\mathbf{x}) = P(Y = +1 | X = \mathbf{x})$, where $X$ and $Y$ are two random variables whose joint distribution is $p(\mathbf{x}, \mathbf{y})$. The optimal Bayes rule for the minimization of the expected misclassification rate is

$$(3.7) \quad BR(\mathbf{x}) = \begin{cases} +1, & \text{if } p_1(\mathbf{x}) > \frac{1}{2}, \\ -1, & \text{if } p_1(\mathbf{x}) < \frac{1}{2}. \end{cases}$$

On one hand, from the previous section we know that the minimization of problem (3.1) guarantees (via regularization theory) that the empirical risk $\frac{1}{n}\sum_{i=1}^n (1 - y_i f(\mathbf{x}_i))_+$ converges to the expected error $E[(1 - Y f(x))_+]$. On the other hand, in [40] it is shown that the solution to the problem $\min_f E[(1 - Y f(x))_+]$ is $f^*(\mathbf{x}) = \text{sign}(p_1(\mathbf{x}) - 1/2)$, an equivalent formulation of (3.7). Therefore, the minimizer sought by SVMs is exactly the Bayes rule.

In [41] it is pointed out that if the smoothing parameter $\mu$ in (3.1) is chosen appropriately and the approximation capacity of the RKHS is large enough, then the solution to the SVM problem (3.2) approaches the Bayes rule as $n \to \infty$. For instance, in the two examples shown in the next subsection, where the linear kernel $K(\mathbf{x}, \mathbf{y}) = \mathbf{x}^T \mathbf{y}$ is used, the associated RKHS (made up of linear functions) is rich enough to solve the classification problems. A richer RKHS should be used for more complex decision surfaces (see [41]), for instance, the one induced by the Gaussian kernel or those induced by high degree polynomial kernels. Regarding the choice of $\mu$, methods to determine it in an appropriate manner have been proposed by Wahba [79, 80, 82].

## 3.2 Illustrating the Performance with Simple Examples

In this first example we consider a two-class separable classification problem, where each class is made up of 1000 data points generated from a bivariate normal distribution $N(\mu_i, I)$, with $\mu_1 = (0, 0)$ and $\mu_2 = (10, 10)$. Our aim is to illustrate the performance of the SVM in a simple example and, in particular, the behavior of the algorithm for different values of the regularization parameter $C$ in problem (3.3). The identity mapping $\Phi(\mathbf{x}) = \mathbf{x}$ is used. Figure 4(a) illustrates the result for $C = 1$ (for $C > 1$, the same result is obtained). There are exactly three support vectors and the optimal margin separating hyperplane obtained by the SVM is $1.05x + 1.00y - 10.4 = 0$. For $C = 0.01$, seven support vectors are obtained [see Figure 4(b)], and the discriminant line is $1.02x + 1.00y - 10.4 = 0$. For $C = 0.00001$, 1776 support vectors are obtained [88.8% of the sample; see Figure 4(c)] and the separating hyperplane is $1.00x + 1.00y - 13.0 = 0$. The three hyperplanes are very similar to the (normal theory) linear discriminant function $1.00x + 1.00y - 10.0 = 0$. Notice that the smaller $C$ is, the larger the number of support vectors. This is due to the fact that, in problem (3.3), $C$ penalizes the value of the $\xi_i$ variables, which determine the width of the band that contains the support vectors.

This second example is quite similar to the previous one, but the samples that correspond to each class are not separable. In this case the mean vectors of the two normal clouds (500 data points in each group) are $\mu_1 = (0, 0)$ and $\mu_2 = (4, 0)$, respectively. The theoretical Bayes error is 2.27%. The normal theory (and optimal) separating hyperplane is $x = 2$, that is, $0.5x + 0y - 1 = 0$. The SVM estimated hyperplane (taking $C = 2$) is $0.497x - 0.001y - 1 = 0$. The error on a test data set with 20,000 data points is 2.3%. Figure 4(d) shows the estimated hyperplane and the support vectors (the black points), which represent 6.3% of the sample. To show the behavior of the method when the parameter $C$ varies, Figure 4(e) shows the separating hyperplanes for 30 SVMs that vary $C$ from 0.01 up to 10. All of them look very similar. Finally, Figure 4(f) shows the same 30 hyperplanes when two outlying points (enhanced in black) are added to the left cloud. Since the estimated SVM discriminant functions depend only on the support vectors, the hyperplanes remain unchanged.

### 3.3 The Waveform Data Set

We next illustrate the performance of SVMs on a well-known three-class classification example considered to be a difficult pattern recognition problem [28], the waveform data set introduced in [13]. For the sake of clarity, we reproduce the data description. Each class is generated from a random convex combination of two of three triangular waveforms, namely, $h_1(i) = \max(6 - |i - 11|, 0)$, $h_2(i) = h_1(i - 4)$ and $h_3(i) = h_1(i + 4)$, sampled at the integers $i \in \{1, \ldots, 21\}$, plus a standard Gaussian noise term. Thus, each data point is represented by $\mathbf{x} = (x_1, \ldots, x_{21})$, where each component is defined by

$$x_i = uh_1(i) + (1 - u)h_2(i) + \varepsilon_i, \quad \text{for Class 1,}$$
$$x_i = uh_1(i) + (1 - u)h_3(i) + \varepsilon_i, \quad \text{for Class 2,}$$
$$x_i = uh_2(i) + (1 - u)h_3(i) + \varepsilon_i, \quad \text{for Class 3,}$$

with $u \sim U(0, 1)$ and $\varepsilon_i \sim N(0, 1)$. A nice picture of sampled waveforms can be found on page 404 of [28]. The waveform data base [available from the UCI repository (data sets available from the University of California, Irvine, at http://kdd.ics.uci.edu/)] contains 5000 instances generated using equal prior probabilities. In this experiment we have used 400 data values for training and 4600 for test. Breiman, Friedman, Olshen and Stone [13] reported a Bayes error rate of 14% for this data set. Since we are handling three groups, we use the "one-against-one" approach, in which $\binom{3}{2}$ binary SVM classifiers are trained and the predicted class is found by a voting scheme: each classifier assigns to each datum a class, being the data point assigned to its most voted class [37]. A first run over ten simulations of the experiment using $C = 1$ in problem (3.3) and the Gaussian kernel $K(x, y) = e^{-\|x - y\|^2/200}$ gave an error rate of 14.6%. To confirm the validity of the result, we have run 1000 replications of the experiment. The average error rate over the 1000 simulations on the training data was 10.87% and the average error rate on the test data was 14.67%. The standard errors of the averages were 0.004 and 0.005, respectively. This result improves any other described in the literature to our knowledge. For instance, the best results described in [28] are provided by FLDA and Fisher FDA (flexible discriminant analysis) with MARS (multivariate adaptive regression splines) as the regression procedure (degree = 1), both achieving a test error rate of 19.1%. Figure 5 shows a principal component analysis (PCA) projection of the waveform data into two dimensions with the misclassified test data points (marked in black) for one of the SVM simulations.

## 4. FURTHER EXAMPLES

In this section we will review some well-known applications of SVMs to real-world problems. In particular, we will focus on text categorization, bioinformatics and image recognition.

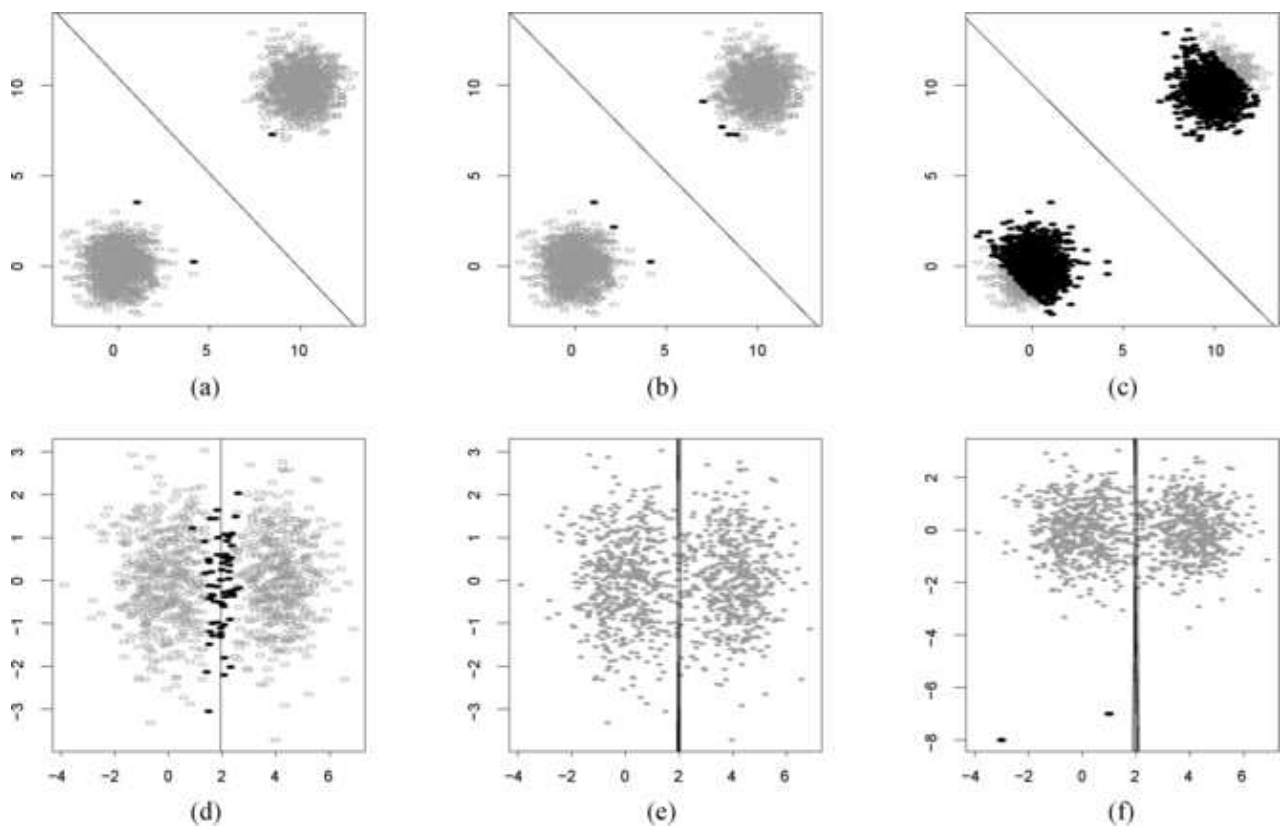

Fig. 4. (a)–(c) *SVM hyperplanes for a separable data set. The support vectors are the black points.* (d)–(f) *SVM hyperplanes for a nonseparable data set.*

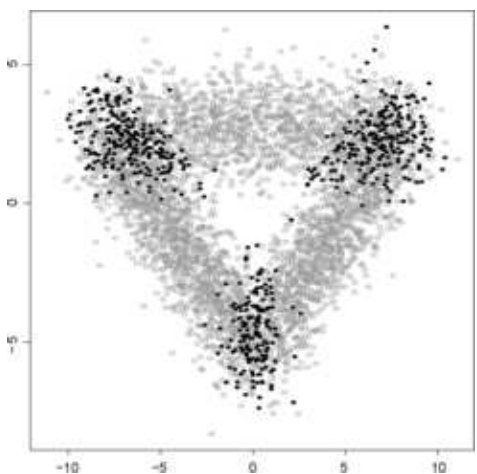

Fig. 5. *A PCA projection of the waveform data. The black points represent the misclassified data points using an SVM with the Gaussian kernel.*

Text categorization consists of the classification of documents into a predefined number of given categories. As an example, consider the document collection made up of Usenet News messages. They are organized in predefined classes such as compu-tation, religion, statistics and so forth. Given a new document, the task is to conduct the category assignment in an automatic way. Text categorization is used by many Internet search engines to select Web pages related to user queries. Documents are represented in a vector space of dimension equal to the number of different words in the vocabulary. Therefore, text categorization problems involve high-dimensional inputs and the data set consists of a sparse document by term matrix. A detailed treatment of SVMs for text categorization can be found in [34]. The performance of SVMs in this task will be illustrated on the Reuters data base. This is a text collection composed of 21,578 documents and 118 categories. The data space in this example has dimension 9947, the number of different words that describe the documents. The results obtained using a SVM with a linear kernel are consistently better along the categories than those obtained with four widely used classification methods: naive Bayes [24], Bayesian networks [29], classification trees [13] and $k$-nearest neighbors [17]. The average rate of success for SVMs is 87% while for the mentioned meth-ods the rates are 72%, 80%, 79% and 82%, respectively (see [34] and [25] for further details). However, the most impressive feature of SVM text classifiers is their training time: SVMs are four times faster than the naive Bayes classifier (the fastest of the other methods) and 35 times faster than classification trees. This performance is due to the fact that SVM algorithms take advantage of sparsity in the document by term matrix. Note that methods that involve the diagonalization of large and dense matrices (like the criterion matrix in FLDA) are out of consideration for text classification because of their expensive computational requirements.

We next outline some SVM applications in bioinformatics. There is an increasing interest in analyzing microarray data, that is, analyzing biological samples using their genetic expression profiles. The SVMs have been applied recently to tissue classification [26], gene function prediction [59], protein subcellular location prediction [31], protein secondary structure prediction [32] and protein fold prediction [23], among other tasks. In almost all cases, SVMs outperformed other classification methods and in their worst case, SVM performance is at least similar to the best non-SVM method. For instance, in protein subcellular location prediction [31], we have to predict protein subcellular positions from prokaryotic sequences. There are three possible location categories: cytoplasmic, periplasmic and extracellular. From a pure classification point of view, the problem reduces to classifying 20-dimensional vectors into three (highly unbalanced) classes. Prediction accuracy for SVMs (with a Gaussian kernel) amounts to 91.4%, while neural networks and a first-order Markov chain [75] have accuracy of 81% and 89.1%, respectively. The results obtained are similar for the other problems. It is important to note that there is still room for improvement.

Regarding image processing, we will overview two well-known problems: handwritten digit identification and face recognition. With respect to the first problem, the U.S. Postal Service data base contains 9298 samples of digits obtained from real-life zip codes (divided into 7291 training samples and 2007 samples for testing). Each digit is represented by a $16 \times 16$ gray level matrix; therefore each data point is represented by a vector in $\mathbb{R}^{256}$. The human classification error for this problem is known to be 2.5% [22]. The error rate for a standard SVM with a third degree polynomial kernel is 4% (see [22] and references therein), while the best known alternative method, the specialized neural network LeNet1 [39], achieves an error rate of 5%. For this problem, using a specialized SVM with a third degree polynomial kernel [22] lowers the error rate to 3.2%—close to the human performance. The key to this specialization lies in the construction of the decision function in three phases: in the first phase, a SVM is trained and the support vectors are obtained; in the second phase, new data points are generated by transforming these support vectors under some groups of transformations, rotations and translations. In the third phase, the final decision hyperplane is built by training a SVM with the new points.

Concerning face recognition, gender detection has been analyzed by Moghaddam and Yang [45]. The data contain 1755 face images (1044 males and 711 females), and the overall error rate for a SVM with a Gaussian kernel is 3.2% (2.1% for males and 4.8% for females). The results for a radial basis neural network [63], a quadratic classifier and FLDA are, respectively, 7.6%, 10.4% and 12.9%.

Another outstanding application of SVMs is the detection of human faces in gray-level images [56]. The problem is to determine in an image the location of human faces and, if there are any, return an encoding of their position. The detection rate for a SVM using a second degree polynomial kernel is 97.1%, while for the best competing system the rate is 94.6%. A number of impressive photographs that show the effectiveness of this application for face location can be consulted in [57].

## 5. EXTENSIONS OF SVMS: SUPPORT VECTOR REGRESSION

It is natural to contemplate how to extend the kernel mapping explained in Section 2 to well-known techniques for data analysis such as principal component analysis, Fisher linear discriminant analysis and cluster analysis. In this section we will describe support vector regression, one of the most popular extensions of support vector methods, and give some references regarding other extensions.

The ideas underlying support vector regression are similar to those within the classification scheme. From an intuitive viewpoint, the data are mapped into a feature space and then a hyperplane is fitted to the mapped data. From a mathematical perspective, the support vector regression function is also derived within the RKHS context. In this case, the loss function involved is known as the $\varepsilon$-insensitive loss function (see [76]), which is defined as $L(y_i, f(\mathbf{x}_i)) = (|f(\mathbf{x}_i) - y_i| - \varepsilon)_+$, $\varepsilon \geq 0$. This loss function ignores errors of size less than $\varepsilon$ (see Figure 6). A discussion of the relationship of the $\varepsilon$-insensitive loss function and the ones used in robust statistics can be found in [28]. Using this loss function, the following optimization problem, similar to (3.1) (also consisting of the minimization of a Tikhonov regularization functional), arises:

$$(5.1) \quad \min_{f \in H_K} \frac{1}{n} \sum_{i=1}^{n} (|f(\mathbf{x}_i) - y_i| - \varepsilon)_+ + \mu \|f\|_K^2,$$

where $\mu > 0$, $H_K$ is the RKHS associated with the kernel $K$, $\|f\|_K$ denotes the norm of $f$ in the RKHS and $(\mathbf{x}_i, y_i)$ are the sample data points.

Once more, by the representer theorem, the solution to problem (5.1) has the form $f(\mathbf{x}) = \sum_{i=1}^{n} \alpha_i \times K(\mathbf{x}_i, \mathbf{x}) + b$, where $\mathbf{x}_i$ are the sample data points. It is immediate to show that $\|f\|_K^2 = \|\mathbf{w}\|^2$, where $\mathbf{w} = \sum_i^n \alpha_i \Phi(\mathbf{x}_i)$ and $\Phi$ is the mapping that defines the kernel function. Thus, problem (5.1) can be restated as

$$(5.2) \quad \min_{\mathbf{w}, b} \frac{1}{n} \sum_{i=1}^{n} (|\mathbf{w}^T \Phi(\mathbf{x}_i) + b - y_i| - \varepsilon)_+ + \mu \|\mathbf{w}\|^2.$$

Since the $\varepsilon$-insensitive loss function is nondifferentiable, this problem has to be formulated so that it can be solved by appropriate optimization methods. Straightforwardly, the equivalent (convex) problem to solve is

$$\min_{\mathbf{w}, b, \xi, \xi'} \quad \frac{1}{2} \|\mathbf{w}\|^2 + C \sum_{i=1}^{n} (\xi_i + \xi_i')$$

$$\text{s.t.} \quad (\mathbf{w}^T \Phi(\mathbf{x}_i) + b) - y_i \leq \varepsilon + \xi_i,$$

$$i = 1, \dots, n,$$

$$(5.3)$$

$$y_i - (\mathbf{w}^T \Phi(\mathbf{x}_i) + b) \leq \varepsilon + \xi_i',$$

$$i = 1, \dots, n,$$

$$\xi_i, \xi_i' \geq 0, \quad i = 1, \dots, n,$$

where $C = 1/(2\mu n)$. Notice that $\varepsilon$ appears only in the constraints, forcing the solution to be calculated by taking into account a confidence band around the regression equation. The $\xi_i$ and $\xi_i'$ are slack variables that allow for some data points to stay outside the confidence band determined by $\varepsilon$. This is the standard support vector regression formulation. Again, the dual of problem (5.3) is a convex quadratic optimization problem, and the regression function takes the same form as equation (2.1). For a detailed exposition of support vector regression, refer to [71] or [69].

One of the most popular applications of support vector regression concerns load forecasting, an important issue in the power industry. In 2001 a proposal based on SVMs for regression was the winner of the European Network of Excellence on Intelligent Technologies competition. The task was to supply the prediction of maximum daily values of electrical loads for January 1999 (31 data values altogether). To this aim each challenger was given half an hour loads, average daily temperatures and the holidays for the period 1997–1998. The mean absolute percentage error for daily data using the SVM regression model was about 2%, significantly improving the results of most competition proposals. It is important to point out that the SVM procedure used in the contest was standard, in the sense that no special modifications were made for the particular problem at hand. See [14] for further details.

Many other kernel methods have been proposed in the literature. To name a few, there are extensions to PCA [70], Fisher discriminant analysis [6, 44], cluster analysis [8, 46], partial least squares [66], time series analysis [50], multivariate density estimation [49, 54, 68], classification with asymmetric proximities [52], combination with neural network models [53] and Bayesian kernel methods [74].

## 6. OPEN ISSUES AND FINAL REMARKS

The underlying model implemented in SVMs is determined by the choice of the kernel. Deciding which kernel is the most suitable for a given application is obviously an important (and open) issue. A possible approach is to impose some restrictions directly on the structure of the classification (or regression) function $f$ implemented by the SVM. A way to proceed is to consider a linear differential operator $D$, and choose $K$ as the Green's function for the operator $D^*D$, where $D^*$ is the adjoint operator of $D$ [4]. It is easy to show that the penalty term $\|f\|_K^2$ equals $\|Df\|_{L_2}^2$. Thus, the choice of the differential operator $D$ imposes smoothing conditions on the solution $f$. This is also the approach used in functional data analysis [65]. For instance, if $D^*D$ is the Laplacian operator, the kernels obtained are harmonic functions. The simplest case corresponds to (see, e.g., [35]) $K(\mathbf{x}, \mathbf{y}) = \mathbf{x}^T \mathbf{y} + \mathbf{c}$, where $c$ is a constant. Another interesting example is the Gaussian kernel. This kernel arises from a differential operator which penalizes an infinite sum of derivatives. The details for its derivation can be found in [63].

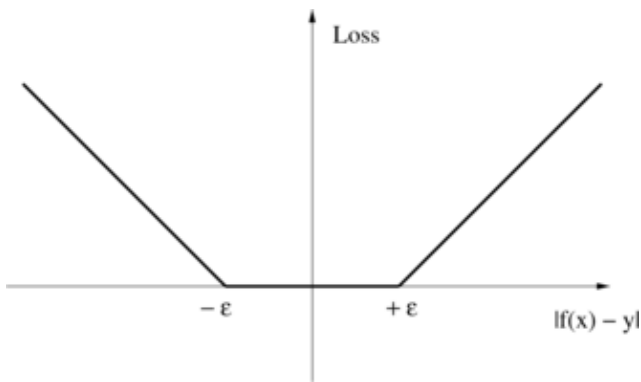

Fig. 6. *The $\varepsilon$-insensitive loss function $L(y_i, f(\mathbf{x}_i)) = (|f(\mathbf{x}_i) - y_i| - \varepsilon)_+$, $\varepsilon > 0$.*

A different approach is to build a specific kernel directly for the data at hand. For instance, Wu and Amari [83] proposed the use of differential geometry methods [2] to derive kernels that improve class separation in classification problems.

An alternative research line arises when a battery of different kernels is available. For instance, when dealing with handwriting recognition, there are a number of different (nonequivalent) metrics that provide complementary information. The task here is to derive a single kernel which combines the most relevant features of each metric to improve the classification performance (see, e.g., [38] or [42]).

Regarding more theoretical questions, Cucker and Smale [21], as already mentioned, provided sufficient conditions for the statistical consistency of SVMs from a functional analysis point of view (refer to the Appendix for the details). On the other hand, the statistical learning theory developed by Vapnik and Chervonenkis (summarized in [77]) provides necessary and sufficient conditions in terms of the Vapnik–Chervonenkis (VC) dimension (a capacity measure for functions). However, the estimation of the VC dimension for SVMs is often not possible and the relationship between both approaches is still an open issue.

From a statistical point of view an important subject remains open: the interpretability of the SVM outputs. Some (practical) proposals can be consulted in [62, 76] and [72] about the transformation of the SVM classification outputs into a posteriori class probabilities.

Regarding the finite sample performance of SVMs, a good starting point can be found in [55], where bias and variability computations for linear inversion algorithms (a particular case of regularization methods) are studied. The way to extend these ideas to the SVM nonlinear case is an interesting open problem.

Concerning software for SVMs, a variety of implementations are freely available from the Web, most reachable at <http://www.kernel-machines.org/>. In particular, Matlab toolboxes and R/Splus libraries can be downloaded from this site. Additional information on implementation details concerning SVMs can be found in [20] and [69].

As a final proposal, a novice reader could find it interesting to review a number of other regularization methods, such as penalized likelihood methods [27], classification and regression with Gaussian processes [72, 82], smoothing splines [81], functional data analysis [65] and kriging [19].

## APPENDIX: STATISTICAL CONSISTENCY OF THE EMPIRICAL RISK

When it is not possible to assume a parametric model for the data, ill-posed problems arise. The number of data points which can be recorded is finite, while the unknown variables are functions which require an infinite number of observations for their exact description. Therefore, finding a solution implies a choice from an infinite collection of alternative models. A problem is well-posed in the sense of Hadamard if (1) a solution exists; (2) the solution is unique; (3) the solution depends continuously on the observed data. A problem is ill-posed if it is not well-posed.

Inverse problems constitute a broad class of ill-posed problems [73]. Classification, regression and density estimation can be regarded as inverse problems. In the general setting, we consider a mapping $H_1 \xrightarrow{A} H_2$, where $H_1$ represents a metric function space and $H_2$ represents a metric space in which the observed data (which could be functions) live. For instance, in a linear regression problem, $H_1$ corresponds to the finite-dimensional vector space $\mathbb{R}^{k+1}$, where $k$ is the number of regressors; $H_2$ is $\mathbb{R}^n$, where $n$ is the number of data points; and $A$ is the linear operator induced by the data matrix of dimension $n \times (k + 1)$. Let $\mathbf{y} = (y_1, \ldots, y_n)$ be the vector of response variables and denote by $f$ the regression equation we are looking for. Then the regression problem consists of solving the inverse problem $Af = \mathbf{y}$. A similar argument applies to the classification setting. In this case, the $\mathbf{y}$ values live in a compact subset of the $H_2$ space [77].

An example of an inverse problem in which $H_2$ is a function space is the density estimation one. In this problem $H_1$ and $H_2$ are both function spaces and $A$ is a linear integral operator given by $(Af)(\mathbf{x}) = \int K(\mathbf{x}, \mathbf{y}) f(\mathbf{y}) \, d\mathbf{y}$, where $K$ is a predetermined kernel function and $f$ is the density function we are seeking. The problem to solve is $Af = F$, where $F$ is the distribution function. If $F$ is unknown, the empirical distribution function $F_n$ is used instead, and the inverse problem to solve is $Af = \mathbf{y}$, with $\mathbf{y} = F_n$.

We will focus on classification and regression tasks. Therefore, we assume there exist a function $f : X \longrightarrow Y$ and a probability measure $p$ defined in $X \times Y$ so that $E[y|\mathbf{x}] = f(\mathbf{x})$. For an observed sample $\{(\mathbf{x}_i, y_i) \in X \times Y\}_{i=1}^{n}$, the goal is to obtain the "best" possible solution to $Af = \mathbf{y}$, where, as mentioned above, $\mathbf{y}$ is the $n$-dimensional vector of $y_i$'s and $A$ is an operator that depends on the $\mathbf{x}_i$ values. To evaluate the quality of a particular solution, a "loss function" $L(f; \mathbf{x}, y)$ has to be introduced, which we will denote $L(y, f(\mathbf{x}))$ in what follows. A common example of a loss function for regression is the quadratic loss $L(y, f(\mathbf{x})) = (y - f(\mathbf{x}))^2$.

Consider the Banach space $C(X)$ of continuous functions on $X$ with the norm $\|f\|_\infty = \sup_{\mathbf{x} \in X} |f(\mathbf{x})|$. The solution to the inverse problem in each case is the minimizer $f^*$ of the risk functional $R(f) : C(X) \longrightarrow \mathbb{R}$ defined by (see [21])

$$(\text{A.1}) \quad R(f) = \int_{X \times Y} L(y, f(\mathbf{x})) p(\mathbf{x}, y) \, d\mathbf{x} \, dy.$$

Of course, the solution depends on the function space in which $f$ lives. Following [21], the hypothesis space, denoted by $\mathcal{H}$ in the sequel, is chosen to be a compact subset of $C(X)$. In particular, only bounded functions $f : X \longrightarrow Y$ are considered.

In these conditions, and assuming a continuous loss function $L$, Cucker and Smale [21] proved that the functional $R(f)$ is continuous. The existence of $f^* = \arg\min_{f \in \mathcal{H}} R(f)$ follows from the compactness of $\mathcal{H}$ and the continuity of $R(f)$. In addition, if $\mathcal{H}$ is convex, $f^*$ will be unique and the problem becomes well-posed.

In practice, it is not possible to calculate $R(f)$ and the empirical risk $R_n(f) = \frac{1}{n} \sum_{i=1}^{n} L(y_i, f(\mathbf{x}_i))$ must be used. This is not a serious complication since asymptotic uniform convergence of $R_n(f)$ to the risk functional $R(f)$ is a proven fact (see [21]).

In summary, imposing compactness on the hypothesis space assures well-posedness of the problem to be solved and uniform convergence of the empirical error to the risk functional for a broad class of loss functions, including the square loss and loss functions used in the SVM setting.

The question of how to impose compactness on the hypothesis space is fixed by regularization theory. A possibility (followed by SVMs) is to minimize Tikhonov's regularization functional

$$(\text{A.2}) \quad \min_{f \in H} \frac{1}{n} \sum_{i=1}^{n} L(y_i, f(\mathbf{x}_i)) + \lambda \Omega(f),$$

where $\lambda > 0$, $H$ is an appropriate function space, and $\Omega(f)$ is a convex positive functional. By standard optimization theory arguments, it can be shown that, for fixed $\lambda$, the inequality $\Omega(f) \leq C$ holds for a constant $C > 0$. Therefore, the space where the solution is searched takes the form $\mathcal{H} = \{f \in H : \Omega(f) \leq C\}$, that is, a convex compact subset of $H$.

## ACKNOWLEDGMENTS

Thanks are extended to Executive Editors George Casella and Edward George, and an anonymous editor for their very helpful comments. The first author was supported in part by Spanish Grants TIC2003-05982-C05-05 (MCyT) and MTM2006-14961-C05-05 (MEC). The second author was supported in part by Spanish Grants SEJ2004-03303 and 06/HSE/0181/2004.

## REFERENCES

[1] Aizerman, M. A., Braverman, E. M. and Rozonoer, L. I. (1964). Theoretical foundations of the potential function method in pattern recognition learning. *Automat. Remote Control* **25** 821–837.

[2] Amari, S.-I. (1985). *Differential-Geometrical Methods in Statistics. Lecture Notes in Statist.* **28**. Springer, New York. MR0788689

[3] Aronszajn, N. (1950). Theory of reproducing kernels. *Trans. Amer. Math. Soc.* **68** 337–404. MR0051437

[4] Aronszajn, N. (1951). Green's functions and reproducing kernels. In *Proc. Symposium on Spectral Theory and Differential Problems* 355–411.

[5] Baeza-Yates, R. and Ribeiro-Neto, B. (1999). *Modern Information Retrieval*. Addison-Wesley, Harlow.

[6] Baudat, G. and Anouar, F. (2000). Generalized discriminant analysis using a kernel approach. *Neural Computation* **12** 2385–2404.

[7] BAZARAA, M. S., SHERALI, H. D. and SHETTY, C. M. (1993). *Nonlinear Programming: Theory and Algorithms*, 2nd ed. Wiley, New York.

[8] BEN-HUR, A., HORN, D., SIEGELMANN, H. and VAPNIK, V. (2001). Support vector clustering. *J. Mach. Learn. Res.* **2** 125–137.

[9] BENNETT, K. P. and CAMPBELL, C. (2000). Support vector machines: Hype or hallelujah? *SIGKDD Explorations* **2** (2) 1–13.

[10] BOSER, B. E., GUYON, I. and VAPNIK, V. (1992). A training algorithm for optimal margin classifiers. In *Proc. Fifth ACM Workshop on Computational Learning Theory* (*COLT*) 144–152. ACM Press, New York.

[11] BOUSQUET, O. and ELISSEEFF, A. (2002). Stability and generalization. *J. Mach. Learn. Res.* **2** 499–526. MR1929416

[12] BREIMAN, L. (2001). Statistical modeling: The two cultures (with discussion). *Statist. Sci.* **16** 199–231. MR1874152

[13] BREIMAN, L., FRIEDMAN, J., OLSHEN, R. and STONE, C. (1984). *Classification and Regression Trees*. Wadsworth, Belmont, CA. MR0726392

[14] CHEN, B.-J., CHANG, M.-W. and LIN, C.-J. (2004). Load forecasting using support vector machines: A study on EUNITE competition 2001. *IEEE Transactions on Power Systems* **19** 1821–1830.

[15] CORTES, C. and VAPNIK, V. (1995). Support-vector networks. *Machine Learning* **20** 273–297.

[16] COVER, T. M. (1965). Geometrical and statistical properties of systems of linear inequalities with applications in pattern recognition. *IEEE Transactions on Electronic Computers* **14** 326–334.

[17] COVER, T. M. and HART, P. E. (1967). Nearest neighbour pattern classification. *IEEE Trans. Inform. Theory* **13** 21–27.

[18] COX, D. and O'SULLIVAN, F. (1990). Asymptotic analysis of penalized likelihood and related estimators. *Ann. Statist.* **18** 1676–1695. MR1074429

[19] CRESSIE, N. (1993). *Statistics for Spatial Data*. Wiley, New York. MR1239641

[20] CRISTIANINI, N. and SHAWE-TAYLOR, J. (2000). *An Introduction to Support Vector Machines*. Cambridge Univ. Press.

[21] CUCKER, F. and SMALE, S. (2002). On the mathematical foundations of learning. *Bull. Amer. Math. Soc.* (*N.S.*) **39** 1–49. MR1864085

[22] DECOSTE, D. and SCHÖLKOPF, B. (2002). Training invariant support vector machines. *Machine Learning* **46** 161–190.

[23] DING, C. and DUBCHAK, I. (2001). Multi-class protein fold recognition using support vector machines and neural networks. *Bioinformatics* **17** 349–358.

[24] DOMINGOS, P. and PAZZANI, M. (1997). On the optimality of the simple Bayesian classifier under zero-one loss. *Machine Learning* **29** 103–130.

[25] DUMAIS, S., PLATT, J., HECKERMAN, D. and SAHAMI, M. (1998). Inductive learning algorithms and representations for text categorization. In *Proc. 7th International Conference on Information and Knowledge Management* 148–155. ACM Press, New York.

[26] FUREY, T. S., CRISTIANINI, N., DUFFY, N., BEDNARSKI, D., SCHUMMER, M. and HAUSSLER, D. (2000). Support vector machine classification and validation of cancer tissue samples using microarray expression data. *Bioinformatics* **16** 906–914.

[27] GREEN, P. J. (1999). Penalized likelihood. *Encyclopedia of Statistical Sciences* **Update 3** 578–586. Wiley, New York.

[28] HASTIE, T., TIBSHIRANI, R. and FRIEDMAN, J. (2001). *The Elements of Statistical Learning*. Springer, New York. MR1851606

[29] HECKERMAN, D., GEIGER, D. and CHICKERING, D. (1995). Learning Bayesian networks: The combination of knowledge and statistical data. *Machine Learning* **20** 197–243.

[30] HERBRICH, R. (2002). *Learning Kernel Classifiers: Theory and Algorithms*. MIT Press, Cambridge, MA.

[31] HUA, S. and SUN, Z. (2001). Support vector machine approach for protein subcellular localization prediction. *Bioinformatics* **17** 721–728.

[32] HUA, S. and SUN, Z. (2001). A novel method of protein secondary structure prediction with high segment overlap measure: Support vector machine approach. *J. Molecular Biology* **308** 397–407.

[33] IVANOV, V. V. (1976). *The Theory of Approximate Methods and their Application to the Numerical Solution of Singular Integral Equations*. Noordhoff International, Leyden. MR0405045

[34] JOACHIMS, T. (2002). *Learning to Classify Text Using Support Vector Machines*. Kluwer, Boston.

[35] KANWAL, R. P. (1983). *Generalized Functions*. Academic Press, Orlando, FL. MR0732788

[36] KIMELDORF, G. S. and WAHBA, G. (1970). A correspondence between Bayesian estimation on stochastic processes and smoothing by splines. *Ann. Math. Statist.* **41** 495–502. MR0254999

[37] KRESSEL, U. (1999). Pairwise classification and support vector machines. In *Advances in Kernel Methods—Support Vector Learning* (B. Schölkopf, C. J. C. Burges and A. J. Smola, eds.) 255–268. MIT Press, Cambridge, MA.

[38] LANCKRIET, G. R. G., CRISTIANINI, N., BARLETT, P., EL GHAOUI, L. and JORDAN, M. I. (2002). Learning the kernel matrix with semi-definite programming. In *Proc. 19th International Conference on Machine Learning* 323–330. Morgan Kaufmann, San Francisco.

[39] LECUN, Y., BOSER, B., DENKER, J. S., HENDERSON, D., HOWARD, R. E., HUBBARD, W. and JACKEL, L. D. (1989). Backpropagation applied to handwritten zip code recognition. *Neural Computation* **1** 541–551.

[40] LIN, Y. (2002). Support vector machines and the Bayes rule in classification. *Data Min. Knowl. Discov.* **6** 259–275. MR1917926

[41] LIN, Y., WAHBA, G., ZHANG, H. and LEE, Y. (2002). Statistical properties and adaptive tuning of sup-

port vector machines. *Machine Learning* **48** 115–136.

[42] MARTIN, I., MOGUERZA, J. M. and MUÑOZ, A. (2004). Combining kernel information for support vector classification. *Multiple Classifier Systems. Lecture Notes in Comput. Sci.* **3077** 102–111. Springer, Berlin.

[43] MERCER, J. (1909). Functions of positive and negative type and their connection with the theory of integral equations. *Philos. Trans. Roy. Soc. London A* **209** 415–446.

[44] MIKA, S., RÄTSCH, G., WESTON, J., SCHÖLKOPF, B. and MÜLLER, K.-R. (1999). Fisher discriminant analysis with kernels. In *Neural Networks for Signal Processing* (Y.-H. Hu, J. Larsen, E. Wilson and S. Douglas, eds.) 41–48. IEEE Press, Piscataway, NJ.

[45] MOGHADDAM, B. and YANG, M.-H. (2002). Learning gender with support faces. *IEEE Trans. Pattern Analysis and Machine Intelligence* **24** 707–711.

[46] MOGUERZA, J. M., MUÑOZ, A. and MARTIN-MERINO, M. (2002). Detecting the number of clusters using a support vector machine approach. *Proc. International Conference on Artificial Neural Networks. Lecture Notes in Comput. Sci.* **2415** 763–768. Springer, Berlin.

[47] MOGUERZA, J. M. and PRIETO, F. J. (2003). An augmented Lagrangian interior-point method using directions of negative curvature. *Math. Program. Ser. A* **95** 573–616. MR1969766

[48] MUKHERJEE, S., RIFKIN, P. and POGGIO, T. (2003). Regression and classification with regularization. *Nonlinear Estimation and Classification. Lecture Notes in Statist.* **171** 111–128. Springer, New York. MR2005786

[49] MUKHERJEE, S. and VAPNIK, V. (1999). Multivariate density estimation: A support vector machine approach. Technical Report, AI Memo 1653, MIT AI Lab.

[50] MÜLLER, K.-R., SMOLA, A. J., RÄTSCH, G., SCHÖLKOPF, B., KOHLMORGEN, J. and VAPNIK, V. (1999). Using support vector machines for time series prediction. In *Advances in Kernel Methods—Support Vector Learning* (B. Schölkopf, C. J. C. Burges and A. J. Smola, eds.) 243–253. MIT Press, Cambridge, MA.

[51] MÜLLER, P. and RIOS INSUA, D. (1998). Issues in Bayesian analysis of neural network models. *Neural Computation* **10** 749–770.

[52] MUÑOZ, A., MARTIN, I. and MOGUERZA, J. M. (2003). Support vector machine classifiers for asymmetric proximities. *Artificial Neural Networks and Neural Information. Lecture Notes in Comput. Sci.* **2714** 217–224. Springer, Berlin.

[53] MUÑOZ, A. and MOGUERZA, J. M. (2003). Combining support vector machines and ARTMAP architectures for natural classification. *Knowledge-Based Intelligent Information and Engineering Systems. Lecture Notes in Artificial Intelligence* **2774** 16–21. Springer, Berlin.

[54] MUÑOZ, A. and MOGUERZA, J. M. (2006). Estimation of high-density regions using one-class neighbor machines. *IEEE Trans. Pattern Analysis and Machine Intelligence* **28** 476–480.

[55] O'SULLIVAN, F. (1986). A statistical perspective on ill-posed inverse problems (with discussion). *Statist. Sci.* **1** 502–527. MR0874480

[56] OSUNA, E., FREUND, R. and GIROSI, F. (1997). Training support vector machines: An application to face detection. In *Proc. IEEE Conference on Computer Vision and Pattern Recognition* 130–136. IEEE Press, New York.

[57] OSUNA, E., FREUND, R. and GIROSI, F. (1997). Support vector machines: Training and applications. CBCL Paper 144/AI Memo 1602, MIT AI Lab.

[58] OSUNA, E., FREUND, R. and GIROSI, F. (1997). An improved training algorithm for support vector machines. In *Proc. IEEE Workshop on Neural Networks for Signal Processing* 276–285. IEEE Press, New York.

[59] PAVLIDIS, P., WESTON, J., CAI, J. and GRUNDY, W. N. (2001). Gene functional classification from heterogeneous data. In *Proc. Fifth Annual International Conference on Computational Biology* 249–255. ACM Press, New York.

[60] PHILLIPS, D. L. (1962). A technique for the numerical solution of certain integral equations of the first kind. *J. Assoc. Comput. Mach.* **9** 84–97. MR0134481

[61] PLATT, J. C. (1999). Fast training of support vector machines using sequential minimal optimization. In *Advances in Kernel Methods—Support Vector Learning* (B. Schölkopf, C. J. C. Burges and A. J. Smola, eds.) 185–208. MIT Press, Cambridge, MA.

[62] PLATT, J. C. (2000). Probabilities for SV machines. In *Advances in Large-Margin Classifiers* (P. J. Bartlett, B. Schölkopf, D. Schuurmans and A. J. Smola, eds.) 61–74. MIT Press, Cambridge, MA.

[63] POGGIO, T. and GIROSI, F. (1990). Networks for approximation and learning. *Proc. IEEE* **78** 1481–1497.

[64] POGGIO, T., MUKHERJEE, S., RIFKIN, R., RAKHLIN, A. and VERRI, A. (2001). *b*. CBCL Paper 198/AI Memo 2001-011, MIT AI Lab.

[65] RAMSAY, J. O. and SILVERMAN, B. W. (1997). *Functional Data Analysis.* Springer, New York.

[66] ROSIPAL, R. and TREJO, L. J. (2001). Kernel partial least squares regression in reproducing kernel Hilbert space. *J. Mach. Learn. Res.* **2** 97–123.

[67] SCHÖLKOPF, B., HERBRICH, R., SMOLA, A. J. and WILLIAMSON, R. C. (2001). A generalized representer theorem. *Lecture Notes in Artificial Intelligence* **2111** 416–426. Springer, Berlin.

[68] SCHÖLKOPF, B., PLATT, J. C., SHAWE-TAYLOR, J., SMOLA, A. J. and WILLIAMSON, R. C. (2001). Estimating the support of a high-dimensional distribution. *Neural Computation* **13** 1443–1471.

[69] SCHÖLKOPF, B. and SMOLA, A. J. (2002). *Learning with Kernels.* MIT Press, Cambridge, MA.

[70] SCHÖLKOPF, B., SMOLA, A. J. and MÜLLER, K.-R. (1999). Kernel principal component analysis. In *Advances in Kernel Methods—Support Vector Learn-*

*ing* (B. Schölkopf, C. J. C. Burges and A. J. Smola, eds.) 327–352. MIT Press, Cambridge, MA.

[71] SMOLA, A. J. and SCHÖLKOPF, B. (1998). A tutorial on support vector regression. NeuroColt2 Technical Report Series, NC2-TR-1998-030.

[72] SOLLICH, P. (2002). Bayesian methods for support vector machines: Evidence and predictive class probabilities. *Machine Learning* **46** 21–52.

[73] TIKHONOV, A. N. and ARSENIN, V. Y. (1977). *Solutions of Ill-Posed Problems.* Wiley, New York.

[74] TIPPING, M. (2001). Sparse Bayesian learning and the relevance vector machine. *J. Mach. Learn. Res.* **1** 211–244. MR1875838

[75] VAN KAMPEN, N. G. (1981). *Stochastic Processes in Physics and Chemistry.* North-Holland, Amsterdam. MR0648937

[76] VAPNIK, V. (1995). *The Nature of Statistical Learning Theory.* Springer, New York. MR1367965

[77] VAPNIK, V. (1998). *Statistical Learning Theory.* Wiley, New York. MR1641250

[78] VAPNIK, V. and CHERVONENKIS, A. (1964). A note on a class of perceptrons. *Automat. Remote Control* **25** 103–109.

[79] WAHBA, G. (1980). Spline bases, regularization, and generalized cross validation for solving approximation problems with large quantities of noisy data. In *Approximation Theory III* (W. Cheney, ed.) 905–912. Academic Press, New York. MR0602818

[80] WAHBA, G. (1985). A comparison of GCV and GML for choosing the smoothing parameter in the generalized spline smoothing problem. *Ann. Statist.* **13** 1378–1402. MR0811498

[81] WAHBA, G. (1990). *Spline Models for Observational Data.* SIAM, Philadelphia. MR1045442

[82] WAHBA, G. (1999). Support vector machines, reproducing kernel Hilbert spaces and the randomized GACV. In *Advances in Kernel Methods—Support Vector Learning* (B. Schölkopf, C. J. C. Burges and A. J. Smola, eds.) 69–88. MIT Press, Cambridge, MA.

[83] WU, S. and AMARI, S.-I. (2002). Conformal transformation of kernel functions: A data-dependent way to improve support vector machine classifiers. *Neural Processing Letters* **15** 59–67.